\documentclass[12pt]{article}
\usepackage{mathrsfs}
\usepackage{amsfonts}
\usepackage{amscd}
\usepackage{fancyhdr}
\usepackage{graphicx}

\usepackage{amssymb,amsmath,latexsym,amscd}
\usepackage{pb-diagram,lamsarrow,pb-lams}
\usepackage[usenames]{color}

\let\goth\mathfrak

\def\gg{\goth g}

\def\gs{\goth s}

\def\gz{\goth z}

\def\L{\mathcal L}

\def\bAut{\textbf{Aut}}

\def\beq{\begin{equation}}
\def\eeq{\end{equation}}
\def\bea{\begin{eqnarray}}
\def\eea{\end{eqnarray}}
\def\beas{\begin{eqnarray*}}
\def\eeas{\end{eqnarray*}}

\def\cplus{\hbox{$\supset${\raise1.05pt\hbox{\kern -0.55em
${\scriptscriptstyle +}$}}\ }}

\DeclareMathOperator{\Hom}{Hom}
\DeclareMathOperator{\Aut}{Aut}

\newtheorem{lemma}[equation]{Lemma}

\newtheorem{corollary}[equation]{Corollary}

\newtheorem{proposition}[equation]{Proposition}

\newtheorem{remark}[equation]{Remark}

\title{Universal Central Extensions of Twisted Forms of Split Simple Lie Algebras over Rings}
\author{Jie Sun$\hbox{\,}$\vspace{0.3cm}\\$\hbox{\,}$ {\small University of Alberta, Department of Mathematical and Statistical Sciences},\\ {\small 632 CAB, University of Alberta, Edmonton, Alberta, Canada T6G 2G1}\\ {\small Email:\ jsun@math.ualberta.ca}}
\date{}

\begin{document}
\maketitle

\begin{small}
\noindent {\bf Abstract.} We give sufficient conditions for the
descent construction to be the universal central extension of a
twisted form of a split simple Lie algebra over a ring. In
particular, the universal central extensions of twisted multiloop
Lie tori are obtained by the descent construction.

\bigskip




\end{small}

\bigskip

\section{Introduction}
Central extensions play a crucial role in physics as they can reduce
the study of projective representations to the study of true
representations.  An important example of this is the Witt and
Virasoro algebras which are infinite dimensional Lie algebras with
many applications to physics. They often appear in problems with
conformal symmetry where the essential spacetime is one or two
dimensional and space is periodic, i.e. compactified to a circle. An
example of such a setting is string theory where the string
worldsheet is two dimensional and cylindrical in the case of closed
strings (see \S 4.3 in \cite{G}). Such worldsheets are Riemann
surfaces which are invariant under conformal transformations. The
algebra of infinitesimal conformal transformations is the direct sum
of two copies of the Witt algebra. The Virasoro algebra is a one
dimensional central extension (in this case, the universal central
extension) of the Witt algebra.

The study of projective representations of the Witt algebra can be
reduced to the study of true representations of the Virasoro
algebra. The representations of the Virasoro algebra that are of
interest in most physical applications are the unitary irreducible
highest weight representations. These are completely characterized
by the central charge and the conformal weight corresponding to the
highest weight vector (see \S 3.2 in \cite{KR}). To each affine
Kac-Moody algebra there is an associated Virasoro algebra by
Sugawara's construction (see \S 3.2.3 in \cite{G}). A given unitary
representation of the Kac-Moody algebra then naturally transforms
into a unitary representation of the associated Virasoro algebra.

Kac's loop construction realizes all affine Kac-Moody algebras as
the universal central extensions of loop algebras based on finite
dimensional simple Lie algebras (\cite{K}). Extended affine Lie
algebras (EALAs), which arose in the work of K. Saito and P. Slodowy
on elliptic singularities and in the paper by the physicists R.
H{\o}egh-Krohn and B. Torresani (\cite{H-KT}) on Lie algebras of
interest to quantum gauge field theory, are natural generalizations
of affine Kac-Moody algebras. A mathematical foundation of the
theory of EALAs is provided in \cite{AABGP}. Kac's loop construction
gives inspiration to the study of EALAs. The centreless cores of
extended affine Lie algebras have been characterized axiomatically
as centreless Lie tori. In \cite{N2} E. Neher realizes all EALAs as
central extensions of centreless Lie tori. Almost all centreless Lie
tori, namely those which are finitely generated over their centroids
(f.g.c. for short), can be realized as multiloop Lie algebras based
on finite dimensional simple Lie algebras (\cite{ABFP1} and
\cite{ABFP2}). Using Grothendieck's descent formalism allows us to
view multiloop Lie algebras as twisted forms (\cite{GP1}, \cite{GP2}
and \cite{P2}). This new perspective presents a beautiful bridge
between infinite dimensional Lie theory and descent theory. In
\cite{PPS} a natural construction for central extensions of twisted
forms of split simple Lie algebras over rings is given by using
Galois descent.

The purpose of this article is to study the universal central
extensions of infinite dimensional Lie algebras. In the affine
Kac-Moody case, the universal central extension is one dimensional.
For the ``higher nullity" EALAs, the universal central extensions
are infinite dimensional (\cite{MRY} and \cite{EF}). In \cite{Ka} C.
Kassel constructs the universal central extensions of untwisted
multiloop Lie algebras by using K\"{a}hler differentials. It is much
more complicated in the twisted case. Kassel's model has been
generalized in \cite{BK} under certain conditions. Unfortunately
twisted multiloop Lie tori do not satisfy these conditions. In
\cite{N2} E. Neher constructs central extensions of centreless Lie
tori by using centroidal derivations and states that the graded dual
of the algebra of skew centroidal derivations gives the universal
central extension of a centreless Lie torus. However, it is
difficult to calculate the centroidal derivations in general. In
this article, we give sufficient conditions for the descent
construction in \cite{PPS} to give the universal central extensions
of twisted forms of split simple Lie algebras over rings. In
particular, the universal central extensions of twisted multiloop
Lie tori are given by the descent construction and a good
understanding of the centre is provided.

Throughout $k$ will denote a field of characteristic $0$, and $\gg$
a finite dimensional split simple Lie algebra over $k$. Let $R$ and
$S$ be commutative, associative, unital $k$-algebras. We write
$\gg_{R}=\gg\otimes_{k} R$ and $\gg_{S}=\gg\otimes_{k} S$.

\section{Descent constructions for central extensions}
In this section we will recall Kassel's construction for the
universal central extension of $\gg_{R}$ and the descent
construction for central extensions of twisted forms of $\gg_{R}$.

Let $\L$ be a Lie algebra over $k$ and $V$ a $k$-space. Any cocycle
$P \in Z^2(\L,V)$, where $V$ is viewed as a trivial $\L$-module,
leads to a central extension

$$0\longrightarrow V\longrightarrow \L_P \buildrel \pi
\over \longrightarrow \L\longrightarrow 0$$ of $\L$ by $V$. As a
space $ \L_P = \L\oplus V$, and the bracket $[\, , \,]_P$ on $\L_P$
is given by

$$[x\oplus u,y\oplus v]_{P}=[x,y]\oplus P(x,y) \,\, \text{for} \,\, x,y\in \L \,\, \text{and} \,\, u,v\in V.$$
The equivalence class of this extension depends only on the class of
$P$ in $H^2(\L,V)$, and this gives a parametrization of all
equivalence classes of central extensions of $\L$ by $V$ (see for
example \cite{MP} or \cite{We} for details). In this situation,  we
will henceforth naturally identify $V$ with a subspace of $\L_P$.
Assume $\L$ is perfect. We fix once and for all a universal central
extension $0 \longrightarrow V \longrightarrow \widehat{\L}
\buildrel \pi \over \longrightarrow \L \longrightarrow 0 $
(henceforth referred to as {\it the} universal central extension of
$\L$). We will find it useful at times to think of this extension as
being given by a (fixed in our discussion)  ``universal" cocycle
$\widehat{P}$, thus $\widehat{\L} = \L_{\widehat{P}} = V \oplus \L$.
This cocycle is of course not unique, but the class of $\widehat{P}$
in $H^2(\L,V)$ is unique.

We view $\gg_R$ as a Lie algebra {\it over} $k$ (in general infinite
dimensional) by means of the unique bracket satisfying
\begin{equation}
[x\otimes a,y\otimes b]=[x,y]\otimes ab
\end{equation}
for all $x,y \in \gg$ and $a,b \in R$. Of course $\gg_R$ is also
naturally an $R$-Lie algebra (which is free of finite rank). It will
be clear at all times which of the two structures is being
considered.

Let $(\Omega_{R/k}, d_R)$ be the  $R$-module of K\"{a}hler
differentials of the $k$-algebra $R$. When no confusion is possible,
we will simply write $(\Omega_R, d)$. Following Kassel \cite{Ka}, we
consider the $k$-subspace $dR$ of $\Omega_R$, and the corresponding
quotient map $^{\overline{\,\,\,}} \, : \, \Omega_R \rightarrow
\Omega_R/dR$. We then have a unique  cocycle $\widehat{P} =
\widehat{P}_R \in Z^2(\gg_R , \Omega_R/dR)$ satisfying
\begin{equation}
\widehat{P}(x\otimes a,y\otimes b)=(x|\, y)\overline{adb},
\end{equation}
where $(\cdot\mid\cdot)$ denotes the Killing form of $\gg$.

Let $\widehat{\gg_R}$ be the unique Lie algebra over $k$ with the
underlying space $\gg_{R}\oplus \Omega_{R}/dR$, and the unique
bracket satisfying
\begin{equation}
[x\otimes a,y\otimes b]_{\widehat{P}} = [x,y]\otimes ab \oplus (x|\,
y)\overline{adb}.
\end{equation}
As the notation suggests,
$$
0 \longrightarrow \Omega_{R}/dR \longrightarrow \widehat{\gg_R}
\buildrel \pi \over \longrightarrow \gg_R  \longrightarrow 0
$$
is the universal central extension of $\gg_R$. There are other
different realizations of the universal central extension (see
\cite{N1}, \cite{MP} and \cite{We} for details on three other
different constructions), but Kassel's model is perfectly suited for
our purposes.

We now turn our attention to twisted forms of $\gg_R$ for the flat
topology of $R$, i.e. we look at $R$-Lie algebras $\L$ for which
there exists a faithfully flat and finitely presented extension
$S/R$ such that
\begin{equation}\label{twistedform}
\L\otimes_{R}S \simeq \gg_R\otimes_{R}S \simeq \gg \otimes_k S,
\end{equation}
where the above are isomorphisms of $S$-Lie algebras.

Let $\bAut(\gg)$ be the $k$-algebraic group of automorphisms of
$\gg$. The $R$-group $\bAut(\gg)_R$ obtained by base change is
clearly isomorphic to $\bAut(\gg_R)$. It is an affine, smooth, and
finitely presented group scheme over $R$ whose functor of points is
given by
\begin{equation}
\bAut(\gg_R)(S) =\Aut_S(\gg_R \otimes_{R}S) \simeq \Aut_S(\gg
\otimes_k S).
\end{equation}
 By Grothendieck's theory of descent (see Chapter I \S2 in \cite{Mln}
 and Chapter XXIV in
\cite{SGA}), we have a natural  bijective map
\begin{equation}
\text{\rm Isomorphism classes of twisted forms of} \,\, \gg_R
\longleftrightarrow H^1_{\acute et}\big(R, \bAut(\gg_R)\big).
\end{equation}

The descent construction for central extensions of twisted forms of
$\gg_{R}$ relies on the following fundamental fact about lifting
automorphisms to central extensions.

\begin{proposition} \label{fundamental}
Let $\L$ be a perfect Lie algebra over $k$. Then

(1) There exists a (unique up to equivalence) universal central
extension
$$
0 \longrightarrow V \longrightarrow \widehat{\L} \buildrel \pi \over
\longrightarrow \L \longrightarrow 0.
$$

(2) If  $\L$ is centreless, the centre $\gz(\widehat{\L})$ of
$\widehat{\L}$ is precisely the kernel $V$ of the projection
homomorphism $\pi : \widehat{\L} \rightarrow \L$ above. Furthermore,
the canonical map $Aut_{k}(\widehat{\L})\rightarrow Aut_{k}(\L)$ is
an isomorphism.

\end{proposition}
{\it Proof.} (1) The existence of an initial object in the category
of central extensions of $\L$ is due to Garland \cite{Grl} \S5
Remark 5.11 and Appendix III. (See also Theorem 1.14 in \cite{N1},
\S1.9 Proposition 2 in \cite{MP} and \S7.9 Theorem 7.9.2 in
\cite{We} for details).

(2) This result goes back to van der Kallen (see \S11 in
\cite{vdK}). Other proofs can be found in \cite{N1} Theorem 2.2 and
in \cite{P1} Proposition 2.2, Proposition 2.3 and Corollary 2.1.
\hfill $\square$

\bigskip
We recall the following important observation of lifting
automorphisms of $\gg_{R}$ to its central extensions in \cite{PPS}
Proposition 3.11.

\begin{proposition}\label{Rlift}
Let $\theta \in \Aut_k(\gg_R)$, and let $\widehat{\theta}$ be the
unique lift of $\theta$ to $\widehat{\gg_R}$ (see Proposition
\ref{fundamental}). If $\theta$ is $R$-linear, then
$\widehat{\theta}$ fixes the centre $\Omega_R/dR$ of
$\widehat{\gg_R}$ pointwise. In particular, every $R$-linear
automorphism of $\gg_R$ lifts to every central extension of $\gg_R$.

\end{proposition}

When $S/R$ is a finite Galois ring extension with Galois group $G$,
the descent data corresponding to $\L$, which a priori is an element
of $\bAut(\gg)(S \otimes_R S),$ can now be thought as being given by
a cocycle $u=(u_{g})_{g\in G} \in Z^{1}\big(G,\Aut_{S}(\gg_ S)\big)$
(usual non-abelian Galois cohomology), where the group $G$ acts on
$\Aut_{S}(\gg_ S) = \Aut_{S}(\gg\otimes_k S)$ via
$^{g}{\theta}=(1\otimes g)\circ\theta\circ(1\otimes g^{-1})$. Then
$$\L\simeq\L_{u}=\{X\in\gg_S:u_{g}{}^{g}X=X\ \text{\rm for all}\ g\in G\}.$$
As above, we let $(\Omega_{S},d)$ be the module of K\"{a}hler
differentials of $S/k$ and let $\widehat{\gg_S}=\gg_{S}\oplus
\Omega_{S}/dS$ be the universal central extension of $\gg_{S}$. The
Galois group $G$ acts naturally both on $\Omega_S$ and on the
quotient $k$-space $\Omega_{S}/dS$, in such a way that
$^{g}({\overline{sdt}})=\overline{^{g}sd^{g}t}$. This leads to an
action of $G$ on $\widehat{\gg_{S}}$ for which
$$^{g}{((x\otimes
s)\oplus z)}=(x\otimes\,^{g}{s})\oplus\,^{g}{z}$$ for all $x\in
\gg$, $s\in S$, $z\in\Omega_{S}/dS$, and $g \in G$. One verifies
immediately that the resulting maps are automorphisms of the $k$-Lie
algebra $\widehat{\gg_{S}}$. We
 henceforth identify $G$ with a subgroup of
$\Aut_{k}(\widehat{\gg_{S}})$, and let $G$ act on
$\Aut_{k}(\widehat{\gg_{S}})$ by conjugation, i.e.,
$^{g}{\theta}=g\theta g^{-1}$. Let $\widehat{u}_{g}$ be the unique
lift of $u_{g}$. We recall the descent construction for central
extensions of twisted forms of $\gg_{R}$ in \cite{PPS} Proposition
4.22.

\begin{proposition}\label{twistedcentral}
Let $u=({u_{g}})_{g \in G}$ be a cocycle in
$Z^{1}\big(G,\Aut_{S}(\gg_S)\big)$. Then

(1) $\widehat{u}=(\widehat{u}_{g})_{g \in G}$ is a cocycle in
$Z^{1}\big(G,\Aut_{k}(\widehat{\gg_S})\big)$.

(2) $\L_{\widehat{u}}=\{x\in\widehat{\gg_S}:\
\widehat{u}_{g}{}^{g}{x}=x\ \, \text{\rm for all}\, g\in G\}$ is a
central extension of the descended algebra $\L_u$ corresponding to
$u$.

(3) There exist canonical isomorphisms $\gz(\L_{\widehat{u}}) \simeq
(\Omega_{S}/dS)^G \simeq \Omega_R/dR$.
\end{proposition}

The following proposition in \cite{PPS} Proposition 4.23 gives
equivalent conditions for
$\L_{\widehat{u}}=\L_{u}\oplus\Omega_{R}/dR$.

\begin{proposition}\label{twistedlift}
With the above notation, the following conditions are equivalent.

(1) $\L_{\widehat{u}}=\L_{u}\oplus\Omega_{R}/dR$ and $\L_{u}$ is
stable under the action of the Galois group $G.$

(2) $\widehat{u}_{g}(\L_{u})\subset\L_{u}$ for all $g\in G$.

If these conditions hold, then every $\theta\in Aut_{R}(\L_{u})$
lifts to an automorphism $\widehat{\theta}$ of $\L_{\widehat{u}}$
that fixes the centre of $\L_{\widehat{u}}$ pointwise.

\end{proposition}

\begin{remark}\label{multiloop}
{\rm Multiloop Lie algebras provide special examples of twisted
forms of $\gg_{R}$ in the sense of Galois descent. Given a finite
dimensional split simple Lie algebra $\gg$ over $k$ and commuting
finite order automorphisms $\sigma_{1},\ldots,\sigma_{n}$ of $\gg$
with $\sigma_{i}^{m_{i}}=1$, the $n$-step multiloop Lie algebra of
$(\gg,\sigma_{1},\ldots,\sigma_{n})$ is defined by
$$L(\gg,\sigma_{1},\ldots,\sigma_{n}):=\bigoplus_{(i_{1}\ldots,i_{n})\in \mathbb{Z}^{n}}\gg_{\overline{i}_{1},\ldots,\overline{i}_{n}}\otimes
t_{1}^{i_{1}/m_{1}}\ldots t_{n}^{i_{n}/m_{n}},$$ where $^{-}:\Bbb
Z\rightarrow\Bbb Z/m_{j}\Bbb Z$ is the canonical map for $1\leq
j\leq n$ and
$$\gg_{\overline{i}_{1},\ldots,\overline{i}_{n}}=\{x\in\gg:
\sigma_{j}(x)=\zeta_{m_{j}}^{i_{j}}x\ \text{\rm for}\ 1\leq j\leq
n\}$$ is the simultaneous eigenspace corresponding to the
eigenvalues $\zeta_{m_{j}}$ (the primitive $m_{j}^{th}$ roots of
unity) for $1\leq j\leq n$. A multiloop Lie algebra
$L(\gg,\sigma_{1},\ldots,\sigma_{n})$ is infinite dimensional over
the given base field $k$, but is finite dimensional over its
centroid $R=k[t_{1}^{\pm1},\ldots,t_{n}^{\pm1}]$. Let $S/R$ be the
finite Galois ring extension with
$S=k[t_{1}^{\pm1/m_{1}},\ldots,t_{n}^{\pm1/m_{n}}]$, then the
following $S$-Lie algebra isomorphism
$$L(\gg,\sigma_{1},\ldots,\sigma_{n})\otimes_{R}S\simeq\gg_{R}\otimes_{R}S$$ tells that
$L(\gg,\sigma_{1},\ldots,\sigma_{n})$ is a twisted form of
$\gg_{R}$. This perspective of viewing multiloop Lie algebras as
twisted forms, which is developed in \cite{GP1}, \cite{GP2} and
\cite{P2}, provides a new way to look at their structure through the
lens of descent theory. Thus a multiloop Lie algebra
$L(\gg,\sigma_{1},\ldots,\sigma_{n})$ as a twisted form of $\gg_{R}$
must be isomorphic to an $R$-Lie algebra $\L_{u}$ for some cocycle
$u$ in $Z^{1}\big(G, \Aut_{S}(\gg_{S})\big)$. From the general
theory about the nature of multiloop Lie algebras as twisted forms
(see \cite{P2} Theorem 2.1 for loop algebras, and \cite{GP2} \S5 for
multiloop algebras in general), the cocycle $u=({u_{g}})_{g \in G}$
is constant (i.e., it has trivial Galois action) with
$u_{g}=v_{g}\otimes id$ for all $g\in G$. The multiloop Lie algebra
$\L_u$ then has a basis consisting of eigenvectors of the $u_g$'s,
and therefore the second equivalent condition of Proposition
\ref{twistedlift} holds. Thus for mulitloop Lie algebras, we have
$\L_{\widehat{u}}=\L_{u}\oplus \Omega_{R}/dR$.}
\end{remark}

\section{Universal Central Extensions}
Since $\L_u$ in the above section is perfect (see \S 5.1 and \S 5.2
of \cite{GP2} for details), it admits a universal central extension
$\widehat{\L_u}$. By Proposition \ref{fundamental}, there exists a
canonical map $\widehat{\L_u} \rightarrow \L_{\widehat{u}}$. In this
section, we give a sufficient condition for
$\L_{\widehat{u}}\simeq\widehat{\L_{u}}$. As an application we show
that if $\L_{u}$ is a multiloop Lie torus, then $\L_{\widehat{u}}$
is the universal central extension of $\L_{u}$.

Throughout this section $S/R$ is a finite Galois ring extension with
Galois group $G$. We identify $R$ with a subring of $S$ and
$\Omega_{R}/dR$ with $(\Omega_{S}/dS)^G$ through a chosen
isomorphism. Let $u=({u_{g}})_{g \in G}\in
Z^{1}\big(G,\Aut_{S}(\gg_S)\big)$ be a constant cocycle with
$u_{g}=v_{g}\otimes id$ for all $g\in G$. Then the descended Lie
algebra corresponding to $u$ is
$$\hspace{-0.85in}\L_{u}=\{X\in\gg_{S}:u_{g}{}^{g}X=X\ \text{\rm for all}\
g\in G\}$$ $$\hspace{1in}=\{\Sigma_{i} x_{i}\otimes
a_{i}\in\gg_S:\Sigma_{i}v_{g}(x_{i})\otimes{}^{g}a_{i}=\Sigma_{i}x_{i}\otimes
a_{i}\ \text{\rm for all}\ g\in G\}.$$ Let
$\gg_{0}=\{x\in\gg:v_{g}(x)=x\ \text{\rm for all}\ g\in G\}$. Then
$\gg_{0}$ is a $k$-Lie subalgebra of $\gg$. We write
${\gg_{0}}_{R}=\gg_{0}\otimes_{k} R$. Clearly ${\gg_{0}}_{R}$ is a
$k$-Lie subalgebra of $\L_{u}$. Assume $\gg_{0}$ is perfect and let
$\widehat{{\gg_{0}}_{R}}={\gg_{0}}_{R}\oplus\Omega_{R}/dR$ be the
universal central extension of ${\gg_{0}}_{R}$.

We first prove a useful lemma and then generalize C. Kassel's proof
in \cite{Ka} that $\widehat{\gg_{R}}$ is the universal central
extension of $\gg_{R}$.

\begin{lemma}\label{useful} Let $\L$ be a Lie algebra over $k$ and
let $V$ be a trivial $\L$-module. If $\gs\subset \L$ is a finite
dimensional semisimple $k$-Lie subalgebra and $\L$ is a locally
finite $\gs$-module, then every cohomology class in $H^{2}(\L,V)$
can be represented by an $\gs$-invariant cocycle.
\end{lemma}
{\it Proof.} For any cocycle $P\in Z^{2}(\L,V)$, our goal is to find
another cocycle $P'\in Z^{2}(\L,V)$ such that $[P]=[P']$ and
$P'(\L,\gs)=\{0\}$. Note that $\Hom_{k}(\L,V)$ is a $\L$-module
given by $y.\beta(x)=\beta(-[y,x])$.

Define a $k$-linear map $f:\gs\rightarrow\Hom_{k}(\L,V)$ by
$f(y)(x)=P(x,y)$. We claim that $f\in Z^{1}\big(\gs,
\Hom_{k}(\L,V)\big)$. Indeed, since $P\in Z^{2}(\L,V)$, we have
$$P(x, y)=-P(y, x)\text{ and } P([x, y], z)+P([y, z], x)+P([z, x], y)=0$$ for all $x,y,z\in\L$.
Then $P(x, [y,z])=P([x,y], z)+P([z,x], y)$, namely
$f([y,z])(x)=f(z)([x,y])+f(y)([z,x])$ for all $x,y,z\in\L$. Thus
$$f([y,z])=y.f(z)-z.f(y)$$
implies $f\in Z^{1}\big(\gs, \Hom_{k}(\L,V)\big)$.

By our assumption that $\gs$ is finite dimensional and semisimple,
the Whitehead's first lemma (see \S7.8 in \cite{We}) yields
$H^{1}\big(\gs, \Hom_{k}(\L,V)\big)=0$. Note that the standard
Whitehead's first lemma holds for finite dimensional $\gs$-modules.
However, $\Hom_{k}(\L,V)$ is a direct sum of finite dimensional
$\gs$-modules when $\L$ is a locally finite $\gs$-module and $V$ is
a trivial $\L$-module, so the result easily extends. So
$f=d^{0}(\tau)$ for some $\tau\in\Hom_{k}(\L,V)$, where $d^{0}$ is
the coboundary map from $\Hom_{k}(\L,V)$ to
$C^{1}\big(\gs,\Hom_{k}(\L,V)\big)$.

Let $P'=P+d^{1}(\tau)$, where $d^{1}$ is the coboundary map from
$\Hom_{k}(\L,V)$ to $C^{2}(\L,V)$. Then $[P']=[P]$. For all $x\in\L$
and $y\in\gs$ we have
\begin{eqnarray*}
P'(x,y)&=&P(x,y)+d^{1}(\tau)(x,y)\\
&=&P(x,y)-\tau([x,y])\\
&=&P(x,y)-f(y)(x)=0.
\end{eqnarray*}
Thus $P'$ is an $\gs$-invariant cocycle.

\hfill $\square$

\begin{proposition}\label{genKassel}
Let $\L_{u}$ be the descended algebra corresponding to a constant
cocycle $u=({u_{g}})_{g \in G}\in Z^{1}\big(G,\Aut_{S}(\gg_S)\big)$.
Let $\L_{P}$ be a central extension of $\L_{u}$ with cocycle $P\in
Z^{2}(\L_{u},V)$. Assume $\gg_{0}$ is central simple, then there
exist a $k$-Lie algebra homomorphism
$\psi:\widehat{{\gg_{0}}_{R}}\rightarrow\L_{p}$ and a $k$-linear map
$\varphi:\Omega_{R}/dR\rightarrow V$ such that the following
commutative diagram.
\[
\begin{diagram}
\node{0}\arrow{e}\node{\Omega_{R}/dR}\arrow{e}\arrow{s,l}{\varphi}\node{\widehat{{\gg_{0}}_{R}}}\arrow{e}\arrow{s,l}{\psi}\node{{\gg_{0}}_{R}}\arrow{s,J}\arrow{e}\node{0}\\
\node{0}\arrow{e}\node{V}\arrow{e}\node{\L_{P}}\arrow{e}\node{\L_{u}}\arrow{e}\node{0}
\end{diagram}
\]
\end{proposition}
{\it Proof.} Our goal is to find $P_{0}\in Z^{2}(\L_{u},V)$ with
$[P_{0}]=[P]$ satisfying
\begin{equation}\label{P0}
P_{0}(x\otimes a, y\otimes 1)=0 \text{\ for all\ } x,y\in\gg_{0}
\text{\ and\ } a\in R.
\end{equation} Applying Lemma \ref{useful} to $\L={\gg_{0}}_{R}$ and
$\gs=\gg_{0}\otimes_{k}k$, it is clear that $\L$ is a locally finite
$\gs$-module and thus we can find an $\gs$-invariant cocycle $P'\in
Z^{2}(\L,V)$, where $P'={P_{\ |}}_{\L\times\L}+d^{1}(\tau)$ for some
$\tau\in\Hom_{k}(\L,V)$. We can extend this $\tau$ to get a
$k$-linear map
$\tau_{0}:\L_{u}={\gg_{0}}_{R}\oplus{\gg_{0}}_{R}^{\bot}\rightarrow
V$ by ${{\tau_{0}}_{\ |}}_{{\gg_{0}}_{R}}=\tau$ and ${{\tau_{0}}_{\
|}}_{{{\gg_{0}}_{R}}^{\bot}}=0$. Let $P_{0}=P+d^{1}(\tau_{0})$,
where $d^{1}$ is the coboundary map from $\Hom_{k}(\L_{u},V)$ to
$C^{2}(\L_{u},V)$. Then $[P_{0}]=[P]$ and it is easy to check that
for all $x,y\in\gg_{0}$ and $a\in R$ we have
\begin{eqnarray*}
P_{0}(x\otimes a,y\otimes 1)&=&P(x\otimes a,y\otimes 1)+d^{1}(\tau_{0})(x\otimes a,y\otimes 1)\\
&=&P(x\otimes a,y\otimes 1)+d^{1}(\tau)(x\otimes a,y\otimes 1)\\
&=&P'(x\otimes a,y\otimes 1)=0.
\end{eqnarray*}

Replace $P$ by $P_{0}$. Since $P\in Z^{2}(\L_{u},V)$, we have
\begin{equation}\label{Pantisym}
P(x\otimes a, y\otimes b)=-P(y\otimes b, x\otimes a),
\end{equation}
\begin{equation}\label{Pantiasso}
P([x\otimes a, y\otimes b], z\otimes c)+P([y\otimes b, z\otimes c],
x\otimes a)+P([z\otimes c, x\otimes a], y\otimes b)=0
\end{equation}
for all $x\otimes a,y\otimes b, z\otimes c\in\L_{u}$. We can define
a $k$-linear map $\Omega_{R}/dR\rightarrow V$ as follows. Fix
$a,b\in R$ and define $\alpha:\gg_{0}\times \gg_{0}\rightarrow V$ by
$\alpha(x,y)=P(x\otimes a, y\otimes b).$ Then with $c=1$ in
(\ref{Pantiasso}) we obtain $P([y,z]\otimes b,x\otimes
a)+P([z,x]\otimes a,y\otimes b)=0$ for all $z\in\gg_{0}$. By
(\ref{Pantisym}) we have
$$P([z,x]\otimes a,y\otimes b)=-P([y,z]\otimes b,x\otimes
a)=P(x\otimes a,[y,z]\otimes b).$$ So
$\alpha([z,x],y)=\alpha(x,[y,z])$. This tells us
$\alpha([x,z],y)=\alpha(x,[z,y])$, namely $\alpha$ is an invariant
bilinear form on $\gg_{0}$. Since $\gg_{0}$ is central simple by our
assumption, $\gg_{0}$ has a unique invariant bilinear form up to
scalars. It follows that there is a unique $z_{a,b}\in V$ such that
for all $x,y\in\gg_{0}$ we have
\begin{equation}\label{Pzab}
P(x\otimes a, y\otimes b)=\alpha(x,y)=(x|y)z_{a,b},
\end{equation}
where $(\cdot\mid\cdot)$ denotes the Killing form of $\gg$. From
(\ref{P0}), (\ref{Pantisym}), (\ref{Pantiasso}) and
$(\cdot\mid\cdot)$ is symmetric we have
\begin{equation}\label{zab}
(i)\ z_{a,1}=0,\\
\ (ii)\ z_{a,b}=-z_{b,a},\\
\ (iii)\ z_{ab,c}+z_{bc,a}+z_{ca,b}=0.
\end{equation}
Then by (ii) and (iii) the map $\varphi:\Omega_{R/k}\simeq
H_{1}(R,R)\simeq R\otimes_{k}R/<ab\otimes c-a\otimes bc+ca\otimes
b>\rightarrow V$ given by $\varphi(adb)=z_{a,b}$ is a well-defined
$k$-linear map. Here $H_{1}$ is the Hochschild homology. By (i)
$\varphi$ induces a well-defined $k$-linear map
$\varphi:\Omega_{R}/dR\rightarrow V$ given by
$\varphi(\overline{adb})=z_{a,b}$.

Finally let $\sigma:\L_{u}\rightarrow\L_{P}$ be any section map
satisfying
\begin{equation}\label{section}
[\sigma(x\otimes a),\sigma(y\otimes b)]_{\L_{P}}-\sigma([x,y]\otimes
ab)=P(x\otimes a, y\otimes b)
\end{equation}
for all $x\otimes a,\ y\otimes b\in\L_{u}$. Define
$\psi:\widehat{{\gg_{0}}_{R}}\rightarrow\L_{P}$ by $\psi(X\oplus
Z)=\sigma(X)\oplus\varphi(Z)$ for all $X\in{\gg_{0}}_{R}$ and
$Z\in\Omega_{R}/dR$. Clearly $\psi$ is a well-defined $k$-linear
map. We claim that $\psi$ is a Lie algebra homomorphism. Indeed, let
$x\otimes a,y\otimes b\in{\gg_{0}}_{R}$, then
$$\psi([x\otimes a,y\otimes b]_{\widehat{{\gg_{0}}_{R}}})=\psi([x,y]\otimes ab\oplus(x|y)\overline{adb})=\sigma([x,y]\otimes ab)+(x|y)z_{a,b},$$
$$[\psi(x\otimes a),\psi(y\otimes b)]_{\L_{P}}=[\sigma(x\otimes a),\sigma(y\otimes b)]_{\L_{P}}=\sigma([x,y]\otimes ab)+P(x\otimes a,y\otimes b).$$
By (\ref{Pzab}) this shows that $\psi$ is a Lie algebra
homomorphism. It is easy to check the following diagram is
commutative.

\[
\begin{diagram}
\node{0}\arrow{e}\node{\Omega_{R}/dR}\arrow{e}\arrow{s,l}{\varphi}\node{\widehat{{\gg_{0}}_{R}}}\arrow{e}\arrow{s,l}{\psi}\node{{\gg_{0}}_{R}}\arrow{s,J}\arrow{e}\node{0}\\
\node{0}\arrow{e}\node{V}\arrow{e}\node{\L_{P}}\arrow{e}\node{\L_{u}}\arrow{e}\node{0}
\end{diagram}
\]


\hfill $\square$

\begin{remark}
{\rm The above proposition generalizes C. Kassel's proof in
\cite{Ka}. When $u$ is a trivial cocycle, we have $\L_{u}=\gg_{R}$
and $\gg_{0}=\gg$. The above proposition shows that
$\widehat{\gg_{R}}$ is the universal central extension of
$\gg_{R}$.}
\end{remark}

To understand the universal central extensions of twisted forms of
$\gg_{R}$, we need to construct a cocycle $P_{0}$ which satisfies a
stronger condition than (\ref{P0}). For each $a\in S\backslash\{0\}$
define $\gg_{a}=\{x\in\gg:v_{g}(x)\otimes{}^{g}a=x\otimes a\
\text{\rm for all}\ g\in G\}.$ Then $\gg_{a}$ is a $k$-subspace of
$\gg$. It is easy to check that $\gg_{a}\subset\gg_{ra}$ for any
$r\in R$ and $\gg_{a}\otimes_{k} Ra$ is a $k$-subspace of $\L_{u}$.

\begin{lemma}\label{ga}
(1) $\gg_{a}=\gg_{0}$ if $a\in R\backslash\{0\}$. In particular,
$\gg_{1}=\gg_{0}$.

(2) $[\gg_{a},\gg_{b}]\subset\gg_{ab}$ for any $a,b\in
S\backslash\{0\}$. $[\gg_{a},\gg_{0}]\subset\gg_{a}$ for any $a\in
S\backslash\{0\}$.

(3) $\Hom_{k}(\gg_{a},V)$ is a $\gg_{0}$-module for any $k$-vector
space $V$, $a\in S\backslash\{0\}$.

\end{lemma}
{\it Proof.} (1) If $a\in R$, then ${}^{g}a=a$ for all $g\in G$.
Thus $\gg_{a}=\{x\in\gg:v_{g}(x)\otimes a=x\otimes a\ \text{\rm for
all}\ g\in G\}.$ Clearly $\gg_{a}\supset\gg_{0}$. On the other hand,
let $x\in\gg_{a}$ and let $\{x_{i}\otimes a_{j}\}_{i\in I,j\in J}$
be a $k$-basis of $\gg_{S}$. Assume $x=\Sigma_{i}\lambda_{i}x_{i}$,
$v_{g}(x)=\Sigma_{i}\lambda_{i}^{g}x_{i}$ and
$a=\Sigma_{i}\mu_{j}a_{j}$. Then $v_{g}(x)\otimes a=x\otimes a$
implies that $\Sigma_{i,j}\lambda_{i}^{g}\mu_{j}(x_{i}\otimes
a_{j})=\Sigma_{i,j}\lambda_{i}\mu_{j}(x_{i}\otimes a_{j})$. Thus
$\lambda_{i}^{g}\mu_{j}=\lambda_{i}\mu_{j}$ for all $i\in I,j\in J$
and $ g\in G$. Since $a\neq 0$, there exists $\mu_{j_{a}}\neq 0$. By
$\lambda_{i}^{g}\mu_{j_{a}}=\lambda_{i}\mu_{j_{a}}$ we get
$\lambda_{i}^{g}=\lambda_{i}$ for all $i\in I$ and $g\in G$. Thus
$x\in\gg_{0}$, so $\gg_{a}=\gg_{0}$.

(2) Let $x\in\gg_{a}$ and $y\in\gg_{b}$. Then
$v_{g}([x,y])\otimes{}^{g}(ab)=[v_{g}(x),v_{g}(y)]\otimes{}^{g}(ab)=[v_{g}(x)\otimes{}^{g}a,v_{g}(y)\otimes{}^{g}b]=[x\otimes
a,y\otimes b]=[x,y]\otimes ab$. Thus $[x,y]\in\gg_{ab}$. For any
$a\in S\backslash\{0\}$ we have
$[\gg_{a},\gg_{0}]=[\gg_{a},\gg_{1}]\subset\gg_{a}$.

(3) Let $y\in\gg_{0}$ and $\beta\in\Hom_{k}(\gg_{a},V)$. Define
$y.\beta(x)=\beta(-[y,x])$. We can check $y.\beta$ is a well-defined
$\gg_{0}$ action.

$\hfill$ $\square$

\begin{proposition}\label{strongerP0}
Let $\L_{u}$ be the descended algebra corresponding to a constant
cocycle $u=({u_{g}})_{g \in G}\in Z^{1}\big(G,\Aut_{S}(\gg_S)\big)$.
Let $\L_{P}$ be a central extension of $\L_{u}$ with cocycle $P\in
Z^{2}(\L_{u},V)$. Assume $\gg_{0}$ is simple, then we can construct
a cocycle $P_{0}\in Z^{2}(\L_{u},V)$ with $[P_{0}]=[P]$ satisfying
$P_{0}(x\otimes a, y\otimes 1)=0 \text{\ for all\ }
x\in\gg_{a},y\in\gg_{0} \text{\ and\ } a\in S$.
\end{proposition}
{\it Proof.} For each $a\in S\backslash\{0\}$, let $\L_{a}$ be the
$k$-Lie subalgebra of $\L_{u}$ generated by the elements in
$(\gg_{a}\otimes_{k}Ra)\cup(\gg_{0}\otimes_{k}k)$. Let
$\gs=\gg_{0}\otimes_{k}k$. By Lemma \ref{ga} (2) we have
$[\gg_{a},\gg_{0}]\subset\gg_{a}$, thus $\L_{a}$ is a locally finite
$\gs$-module. Applying Lemma \ref{useful} to $\L=\L_{a}$ and
$\gs=\gg_{0}\otimes_{k}k$, we can find an $\gs$-invariant cocycle
$P'_{a}\in Z^{2}(\L_{a},V)$, where $P'_{a}={P_{\
|}}_{\L_{a}\times\L_{a}}+d^{1}(\tau_{a})$ for some
$\tau_{a}\in\Hom_{k}(\L_{a},V)$. Let $\{x_{i}\otimes a_{j}\}_{i\in
I,j\in J}$ be a $k$-basis of $\L_{u}$. For each $a_{j}$ choose one
$\tau_{a_{j}}\in\Hom_{k}(\L_{a_{j}},V)$. Note that $x_{i}\otimes
a_{j}\in\L_{u}$ implies $x_{i}\in\gg_{a_{j}}$, thus $x_{i}\otimes
a_{j}\in\L_{a_{j}}$. Define $\tau:\L_{u}\rightarrow V$ to be the
unique linear map such that $\tau(x_{i}\otimes
a_{j})=\tau_{a_{j}}(x_{i})$.

Let $P_{0}=P+d^{1}(\tau)$, where $d^{1}$ is the coboundary map from
$\Hom_{k}(\L_{u},V)$ to $C^{2}(\L_{u},V)$. Then $[P_{0}]=[P]$. For
each $a_{j}\ (j\in J)$ it is easy to check that for all
$x\in\gg_{a_{j}},y\in\gg_{0}$ we have
\begin{eqnarray*}
P_{0}(x\otimes a_{j},y\otimes 1)&=&P(x\otimes a_{j},y\otimes 1)+d^{1}(\tau)(x\otimes a_{j},y\otimes 1)\\
&=&P(x\otimes a_{j},y\otimes 1)+d^{1}(\tau_{a_{j}})(x\otimes a_{j},y\otimes 1)\\
&=&P'_{a_{j}}(x\otimes a_{j},y\otimes 1)=0.
\end{eqnarray*}
Note that our proof does not depend on the choice of $\tau_{a_{j}}$
because $ker(d^{0})=\Hom_{k}(\L_{a_{j}},V)^{\gs}$ and different
choices of $\tau_{a_{j}}$ become the same when restricted to
$[\L_{a_{j}},\gs]$. Thus for any $x\otimes
a=\Sigma_{i,j}\lambda_{ij}x_{i}\otimes a_{j}\in\L_{u}$, we have
$P_{0}(x\otimes a,y\otimes 1)=\Sigma_{ij}P_{0}(x_{i}\otimes
a_{j},y\otimes 1)=0$.

$\hfill$ $\square$

We have the following important observation when $\gg$ has a basis
consisting of simultaneous eigenvectors of $\{v_{g}\}_{g\in G}$.
\begin{lemma}\label{SandR}
Let $\mathscr{B}=\{x_{i}\otimes a_{j}\}_{i\in I,j\in J}$ be a
$k$-basis of $\L_{u}$ with $\{x_{i}\}_{i\in I}$ consisting of
simultaneous eigenvectors of $\{v_{g}\}_{g\in G}$. Take
$x_{i}\otimes a_{j}, x_{l}\otimes a_{k}\in\L_{u}$. If
$0\neq\overline{a_{j}da_{k}}\in\Omega_{R}/dR$, then $a_{j}a_{k}\in
R$ and $[x_{i},x_{l}]\in\gg_{0}$.
\end{lemma}
{\it Proof.} Let $v_{g}(x_{i})=\lambda_{g}^{i}x_{i}$, where
$\lambda_{g}^{i}\in k$. If $x_{i}\otimes a_{j}\in\L_{u}$, then
$x_{i}\in\gg_{a_{j}}$. So
$v_{g}(x_{i})\otimes{}^{g}a_{j}=\lambda_{g}^{i}x_{i}\otimes{}^{g}a_{j}=x_{i}\otimes
a_{j}$. Thus $x_{i}\otimes{}^{g}a_{j}=x_{i}\otimes
(\lambda_{g}^{i})^{-1}a_{j}$, and therefore
$x_{i}\otimes{}(^{g}a_{j}-(\lambda_{g}^{i})^{-1}a_{j})=0$. Since
$x_{i}\neq 0$, we have $^{g}a_{j}-(\lambda_{g}^{i})^{-1}a_{j}=0$,
thus ${}^{g}a_{j}=(\lambda_{g}^{i})^{-1}a_{j}$. Similarly, we can
show that ${}^{g}a_{k}=(\lambda_{g}^{l})^{-1}a_{k}$. So if
$\overline{a_{j}da_{k}}\in\Omega_{R}/dR$, then
$\overline{{}^{g}a_{j}d{}^{g}a_{k}}=\overline{a_{j}da_{k}}$ for all
$g\in G$. Note that
$$
\overline{a_{j}da_{k}}=\overline{{}^{g}a_{j}d{}^{g}a_{k}}=\overline{(\lambda_{g}^{i})^{-1}a_{j}d(\lambda_{g}^{l})^{-1}a_{k}}=(\lambda_{g}^{i})^{-1}(\lambda_{g}^{l})^{-1}\overline{a_{j}da_{k}}=(\lambda_{g}^{i}\lambda_{g}^{l})^{-1}\overline{a_{j}da_{k}}.
$$
So if $\overline{a_{j}da_{k}}\neq0$, then
$(\lambda_{g}^{i}\lambda_{g}^{l})^{-1}=\lambda_{g}^{i}\lambda_{g}^{l}=1$.
Thus
${}^{g}(a_{j}a_{k})={}^{g}a_{j}{}^{g}a_{k}=(\lambda_{g}^{i})^{-1}a_{j}(\lambda_{g}^{l})^{-1}a_{k}=a_{j}a_{k}$.
So $a_{j}a_{k}\in R$ and
$[x_{i},x_{l}]\in[\gg_{a_{j}},\gg_{a_{k}}]\subset\gg_{a_{j}a_{k}}=\gg_{0}$
by Lemma \ref{ga}.

$\hfill$ $\square$

Now we are ready to prove the main result of this section.
\begin{proposition}\label{main}
Let $u=({u_{g}})_{g \in G}\in Z^{1}\big(G,\Aut_{S}(\gg_S)\big)$ be a
constant cocycle with $u_{g}=v_{g}\otimes id$. Let $\L_{u}$ be the
descended algebra corresponding to $u$ and let $\L_{\widehat{u}}$ be
the central extension of $\L_{u}$ obtained by the descent
construction (see Proposition \ref{twistedcentral}). Assume
$\gg_{0}$ is central simple and $\gg$ has a basis consisting of
simultaneous eigenvectors of $\{v_{g}\}_{g\in G}$. Assume
$\L_{\widehat{u}}=\L_{u}\oplus\Omega_{R}/dR$, then
$\L_{\widehat{u}}$ is the universal central extension of $\L_{u}$.
\end{proposition}
{\it Proof.} First of all, $\L_{\widehat{u}}$ is perfect. Indeed,
let $X\oplus Z\in\L_{\widehat{u}}$, where $X\in\L_{u}$ and
$Z\in\Omega_{R}/dR$. Since $\L_{u}$ is perfect, we have
$X=\Sigma_{i}[X_{i},Y_{i}]_{\L_{u}}$ for some
$X_{i},Y_{i}\in\L_{u}$. By the assumption
$\L_{\widehat{u}}=\L_{u}\oplus\Omega_{R}/dR$ we have
$\L_{u}\subset\L_{\widehat{u}}$, then
$X_{i},Y_{i}\in\L_{\widehat{u}}$. Thus
$\Sigma_{i}[X_{i},Y_{i}]_{\L_{\widehat{u}}}=\Sigma_{i}[X_{i},Y_{i}]_{\L_{u}}\oplus
W$ for some $W\in\Omega_{R}/dR$. So $X\oplus
Z=\Sigma_{i}[X_{i},Y_{i}]_{\L_{\widehat{u}}}\oplus(Z-W)$, where
$Z-W\in\Omega_{R}/dR\subset[{\gg_{0}}_{R},{\gg_{0}}_{R}]_{\L_{\widehat{u}}}\subset[\L_{\widehat{u}},\L_{\widehat{u}}]_{\L_{\widehat{u}}}$.
Thus $\L_{\widehat{u}}$ is perfect.

Let $\L_{P}$ be a central extension of $\L_{u}$ with cocycle $P\in
Z^{2}(\L_{u},V)$. By Proposition \ref{strongerP0}, we can assume
that $P(x\otimes a, y\otimes 1)=0 \text{\ for all\ }
x\in\gg_{a},y\in\gg_{0} \text{\ and\ } a\in S$. Let
$\sigma:\L_{u}\rightarrow\L_{P}$ be any section of
$\L_{P}\rightarrow \L_{u}$ satisfying
\begin{equation}\label{section}
[\sigma(x\otimes a),\sigma(y\otimes b)]_{\L_{P}}-\sigma([x,y]\otimes
ab)=P(x\otimes a, y\otimes b)
\end{equation}
for all $x\otimes a,\ y\otimes b\in\L_{u}$. Define
$\psi:\L_{\widehat{u}}\rightarrow\L_{P}$ by $\psi(X\oplus
Z)=\sigma(X)+\varphi(Z)$ for all $X\in\L_{u}$ and
$Z\in\Omega_{R}/dR$, where $\varphi:\Omega_{R}/dR\rightarrow V$ is
the map given by $\varphi(\overline{adb})=z_{a,b}$ in Proposition
\ref{genKassel}. Clearly $\psi$ is a well-defined $k$-linear map. We
claim that $\psi$ is a Lie algebra homomorphism. Indeed, let
$x\otimes a,y\otimes b\in\L_{\widehat{u}}$, then
$$\psi([x\otimes a,y\otimes b]_{\L_{\widehat{u}}})=\psi([x,y]\otimes ab\oplus(x|y)\overline{adb})=\sigma([x,y]\otimes ab)+(x|y)\varphi(\overline{adb}),$$
$$[\psi(x\otimes a),\psi(y\otimes b)]_{\L_{P}}=[\sigma(x\otimes a),\sigma(y\otimes b)]_{\L_{P}}=\sigma([x,y]\otimes ab)+P(x\otimes a,y\otimes b).$$
By (\ref{Pzab}) we have $P(x\otimes a,y\otimes b)=(x|y)z_{a,b}$ for
all $x,y\in\gg_{0}$ and $a,b\in R$. If $a,b\in S\backslash R$, we
have two cases. Since $\psi$ is well-defined, we only need to
consider basis elements in $\L_{u}$. Let $\mathscr{B}=\{x_{i}\otimes
a_{j}\}_{i\in I,j\in J}$ be a $k$-basis of $\L_{u}$ with
$\{x_{i}\}_{i\in I}$ consisting of eigenvectors of the $v_{g}$'s.
Take $x_{i}\otimes a_{j}, x_{l}\otimes a_{k}\in\L_{u}$. If
$0\neq\overline{a_{j}da_{k}}\in\Omega_{R}/dR$, then $a_{j}a_{k}\in
R$ and $[x_{i},x_{l}]\in\gg_{0}$ by Lemma \ref{SandR}. Thus
$[x_{i}\otimes a_{j},x_{l}\otimes
a_{k}]_{\L_{\widehat{u}}}\subset{\gg_{0}}_{R}\oplus\Omega_{R}/dR$.
By Proposition \ref{genKassel} $\psi$ is a Lie algebra homomorphism
in this case. If $0=\overline{a_{j}da_{k}}\in\Omega_{R}/dR$, then
$[x_{i}\otimes a_{j},x_{l}\otimes
a_{k}]_{\L_{\widehat{u}}}=[x_{i}\otimes a_{j}a_{k},x_{l}\otimes
1]_{\L_{\widehat{u}}}$. By Proposition \ref{strongerP0} we have
$P(x_{i}\otimes a_{j}a_{k},x_{l}\otimes 1)=0$. So $\psi$ is a Lie
algebra homomorphism as well in this case. It is easy to check the
following diagram is commutative.

\[
\begin{diagram}
\node{0}\arrow{e}\node{\Omega_{R}/dR}\arrow{e}\arrow{s,l}{\varphi}\node{\L_{\widehat{u}}}\arrow{e}\arrow{s,l}{\psi}\node{\L_{u}}\arrow{s,=}\arrow{e}\node{0}\\
\node{0}\arrow{e}\node{V}\arrow{e}\node{\L_{P}}\arrow{e}\node{\L_{u}}\arrow{e}\node{0}
\end{diagram}
\]

$\hfill$ $\square$

\begin{corollary}\label{uceofmultiloop}
If $\L_{u}$ is a multiloop Lie torus over an algebraically closed
field of characteristic $0$, then $\L_{\widehat{u}}$ is the
universal central extension of $\L_{u}$ and the centre of
$\L_{\widehat{u}}$ is $\Omega_{R}/dR$.
\end{corollary}
{\it Proof.} If $\L_{u}$ is a multiloop Lie torus, by Remark
\ref{multiloop} we have
$\L_{\widehat{u}}=\L_{u}\oplus\Omega_{R}/dR$. By the definition of
multiloop Lie algebras, $\{v_{g}\}_{g\in G}$ is a set of commuting
finite order automorphisms of $\gg$, thus $\gg$ has a basis
consisting of simultaneous eigenvectors of $\{v_{g}\}_{g\in G}$. By
our assumption $\L_{u}$ is a multiloop Lie torus over an
algebraically closed field of characteristic $0$, then we have
 $\gg_{0}$ is central simple (see \S3.2 and \S3.3 in \cite{ABFP2} for details). Thus for mulitloop Lie torus
$\L_{u}$, our construction $\L_{\widehat{u}}$ gives the universal
central extension by Proposition \ref{main} and the centre of
$\L_{\widehat{u}}$ is $\Omega_{R}/dR$ by Proposition
\ref{twistedcentral}. $\hfill$ $\square$

\begin{remark}
{\rm Proposition \ref{main} provides a good understanding of the
universal central extensions of twisted forms corresponding to
constant cocycles. The assumption that $\gg_{0}$ is central simple
is crucial for our proof. As an important application, Corollary
\ref{uceofmultiloop} provides a good understanding of the universal
central extensions of twisted multiloop Lie tori. Recently E. Neher
calculated the universal central extensions of twisted multiloop Lie
tori by using a result on a particular explicit description of the
algebra of derivations of multiloop Lie algebras in \cite{P3}.
Discovering more general conditions under which the descent
construction gives the universal central extension remains an open
problem.}
\end{remark}

\noindent{\bf Acknowledgement}. The author would like to thank the
referee for his/her invaluable comments. Lemma \ref{useful} is
suggested to the author by the referee to simplify the original
proof of Proposition \ref{genKassel} and Proposition
\ref{strongerP0}. The author is also grateful to Professor Arturo
Pianzola for his encouragement and advice.

\end{document}